\documentclass[11pt,a4paper]{article}
\usepackage{pifont}
\usepackage{CJK,epsfig,amsfonts,amsgen,amsmath,amstext,amsbsy,amsopn,amsthm
}
\usepackage{mathrsfs}
\usepackage{float}
\usepackage{epsfig,latexsym,amsfonts,amsmath,amssymb}
\usepackage{color}
\usepackage{calc}
\usepackage{tikz}
\usepackage{verbatim}

\usepackage{comment}
\usepackage{geometry}
\usepackage{color}
\usepackage{calc}
\usepackage{tikz}
\usetikzlibrary{decorations.markings}
\tikzstyle{vertex}=[circle, draw, inner sep=0pt, minimum size=6pt]

\newtheorem{theorem}[subsection]{Theorem}

\newtheorem{remark}[subsection]{Remark}

\newtheorem{claim}[subsection]{Claim}

\newtheorem{fact}[subsection]{Fact}
\usepackage{indentfirst} 
\newcounter{Hcase}

\def\hf{\mathcal{F}}
\def\hg{\mathcal{G}}
\def\hh{\mathcal{H}}
\def\hp{\mathcal{P}}
\def\hs{\mathcal{S}}
\def\hb{\mathcal{B}}
\def\ha{\mathcal{A}}
\def \hht{\mathcal{T}}

\begin{document}
\title{A note on the maximum diversity of intersecting families in the symmetric group}
\author{\small Jian Wang${}^{1}$\thanks{E-mail address: wangjian01@tyut.edu.cn}~ and Jimeng Xiao${}^{2}${\thanks{E-mail address: xiaojm@sustech.edu.cn}}
 \\[2mm]
\small ${}^1$\small Department of Mathematics,\\
\small Taiyuan University of Technology, 
Taiyuan, 030024, P. R. ~China\\[3pt]
\small ${}^2$  Department of Mathematics,\\[-0.8ex]
\small Southern University of Science and Technology, Shenzhen,  518055, P. R. ~China\\}

\baselineskip 20pt

\date{}
\maketitle
 \vspace{4mm}

\begin{abstract}
Let $\mathcal{S}_n$ be the symmetric group on the set $[n]:=\{1,2,\ldots,n\}$. A family $\mathcal{F}\subset \mathcal{S}_n$ is called intersecting if  for every $\sigma,\pi\in \mathcal{F}$ there exists some $i\in [n]$ such that $\sigma(i)=\pi(i)$. Deza and Frankl proved that the largest intersecting family of permutations is the full star, that is, the collection of all permutations with a fixed position. The diversity of an intersecting family $\mathcal{F}$ is defined as the minimum number of permutations in $\mathcal{F}$, whose deletion results in a star. In the present paper, by applying the spread approximation method  developed recently by Kupavskii and Zakharov, we prove that for $n\geq 500$ the diversity of an intersecting subfamily of $\mathcal{S}_n$ is at most $(n-3)(n-3)!$, which is best possible.

\vskip 0.1in \noindent%
 \\[7pt]
\end{abstract}

\newcommand{\lr}[1]{\langle #1\rangle}

\parindent 17pt
\baselineskip 17pt

\section{Introduction}

Let $[n]$ denote the standard $n$-element set $\{1,2,\ldots,n\}$. Let $2^{[n]}$ denote the power set of $[n]$ and $\binom{[n]}{k}$ denote the collection of all $k$-element subsets of $[n]$. A subfamily $\mathcal{F}\subset \binom{[n]}{k}$ is called a {\it $k$-uniform family} or a {\it $k$-graph}. We say that $\mathcal{F}\subset \binom{[n]}{k}$ is   {\it intersecting}  if $F\cap F'\neq \emptyset$ for all $F,F'\in \mathcal{F}$.

One of the most important results in extremal set theory is the  Erd\H{o}s-Ko-Rado Theorem.

\begin{theorem}[\cite{ekr}]
    Let $\hf\subset \binom{[n]}{k}$ be an intersecting family with $n\geq 2k$. Then 
    \[
    |\hf|\leq \binom{n-1}{k-1}. 
    \]
\end{theorem}

For $\mathcal{F}\subset \binom{[n]}{k}$ and $i\in [n]$, let
\[
\mathcal{F}(i) = \left\{F\setminus \{i\}\colon i\in F\in \mathcal{F}\right\},\ \mathcal{F}(\bar{i}) = \left\{F \colon i\notin F\in \mathcal{F}\right\}.
\]
Define the {\it diversity} of a family $\mathcal{F}$ as
\[
\gamma(\mathcal{F}) =\min_{1\leq i\leq n} |\mathcal{F}(\bar{i})|.
\]
Lemons and Palmer \cite{LP} (answering a question of Katona) proved that for $n>n_0(k)$, an intersecting family $\hf\subset \binom{[n]}{k}$ has diversity at most $\binom{n-3}{k-2}$.
This was subsequently improved  by \cite{F17}, \cite{Ku1} and \cite{F2020}. The current record is the following.

\begin{theorem}[\cite{FW2022}]\label{thm-fw-2}
Let $n>36k$. Suppose that $\hf\subset \binom{[n]}{k}$ is intersecting. Then
\begin{align}
\gamma(\hf) \leq \binom{n-3}{k-2}.
\end{align}
\end{theorem} 

In the present paper, we mainly consider the corresponding problem for intersecting families of permutations.  

Let $\mathcal{S}_n$ be the symmetric group on the set $[n]$. A family $\mathcal{F}\subset \mathcal{S}_n$ is called {\it intersecting} if for every  $\sigma,\pi\in \mathcal{F}$ there exists $i\in [n]$ such that $\sigma(i)=\pi(i)$.  

Deza and Frankl proved the Erd\H{o}s-Ko-Rado Theorem for symmetric groups.

\begin{theorem}[\cite{DF}]
Let $\mathcal{F}\subset \mathcal{S}_n$ be an intersecting family. Then 
\[
|\hf| \leq (n-1)!.
\]
\end{theorem}

For $i=1,2,\ldots,n$, let $X_i= \{x_{i1},x_{i2},\ldots,x_{in}\}$. 
A family $\mathcal{F}\subset \mathcal{S}_n$ can be viewed as an $n$-partite $n$-graph on partite sets $X_1,X_2,\ldots,X_n$ with edge set
$\{\{x_{1\sigma(1)},x_{2\sigma(2)},\ldots,x_{n\sigma(n)}\}\colon \sigma \in \mathcal{F}\}$.
From this point of view, a family $\mathcal{F}\subset \mathcal{S}_n$ is also a family of $n$-subsets of $[n^2]$. Thus, one can consider the diversity of intersecting families of symmetric groups as the same as that of ordinary set systems.

By applying the spread approximation method developed recently by Kupavskii and Zakharov, we prove the following result. 

\begin{theorem}\label{main_thm}
Let $\mathcal{F}\subset \mathcal{S}_n$ be an intersecting family. If $n \ge 500$, then
\[
\gamma(\mathcal{F}) \le (n-3)(n-3)!.
\]
\end{theorem}

\begin{remark}
Define the {\it triangle family} $\mathcal{T}(n)$ as
\[
\mathcal{T}(n)=\{F \in \mathcal{S}_n \colon |F \cap \{x_{11},x_{22},x_{33}\}| \geq 2\}.
\]
It  is easy to see that $\mathcal{T}(n)$ is intersecting and $\gamma(\mathcal{T}(n))=(n-3)(n-3)!$, showing that Theorem \ref{main_thm} is best possible.
\end{remark}

For $\mathcal{F}\subset 2^{[n]}$ and $A \subset [n]$, define
\[
\mathcal{F}(A) = \left\{F\setminus A\colon A\subset  F\in \mathcal{F}\right\}.
\]
We say that $\mathcal{F}\subset 2^{[n]}$ is {\it $r$-spread} if $|\mathcal{F}(A)| \leq r^{-|A|} |\mathcal{F}|$ for every $A\subset [n]$.

In \cite{alwz}, Alweiss, Lovett, Wu and Zhang found a very useful result which they applied to prove a breakthrough result concerning the forbidden sunflower problem of  Erd\H{o}s and Rado \cite{er}. Using a variant of this statement due to Tao \cite{tao}, Kupavskii and  Zakharov \cite{kup} developed the spread approximation method, which is a powerful tool in proving intersection theorems \cite{ku1,ku2,ku3}. 
For more intersecting theorems on set systems and permutation groups, see \cite{DF2,difr,elki,elke,ekr,F78,F87,F87-2, keli,kelima,kelish,wz,W84} and the references therein.

We say $W$ is a $p$-random subset of $[n]$ if each element $i\in [n]$ is included in $W$ independently with probability $p$. The following theorem gives a key property of spread families.

\begin{theorem}[\cite{alwz},  \cite{tao}, \cite{sto}]\label{thm-key}
Let $\mathcal{F}\subset \binom{[n]}{\leq k}$  be an $r$-spread family and let $W$ be an $(m\delta)$-random subset of $[n]$. Then
\[
Pr\left[\exists F\subset W \mbox{ for some }F\in \mathcal{F}\right]  \geq 1-\left(\frac{1+H(\delta)}{\log_2(r\delta)}\right)^m k,
\]
where $H(\delta) = -\delta\log_2 \delta -(1-\delta)\log_2(1-\delta)$.
\end{theorem}

The ingredients of the proof of our main theorem include the spread approximation  method introduced by Kupavskii and  Zakharov and a theorem of F\"{u}redi about pseudo sunflowers (see section 3).
The rest of the paper is organized as follows. In the next section, we present some preliminary results. In Section 3, we prove the main theorem.

\section{Preliminary results}

Let $X=X_1\cup X_2\cup \ldots\cup X_n$. For $S\subset X$, we say $S$ is a proper subset of $X$ if $|S\cap X_i|\leq 1$ for all $i=1,2,\ldots,n$. Let $\hp(X)$ be the collection of all  proper subsets of $X$. Let $\hp(X,s)$ be the collection of all  proper subsets of $X$ with size at most $s$. We say $S$ is $r$-spread in $\hf$ if $|\mathcal{F}(S)| \leq r^{-|S|} |\mathcal{F}|$.
Define
\[
\hs(\hf,r) =  \left\{S\in \hp(X)\colon S \mbox{ is not $r$-spread in }\hf \right\}.
\]
and let $\hs^*(\hf,r)$ be the collection of all maximal members of $\hs(\hf,r)$. Now, we show that for intersecting families $\hf$ and large enough $r$, there exists a non-empty proper subset in $\hs^*(\hf,r)$.

\begin{fact}\label{fact-2.1}
Let $\hf \subset \binom{[n]}{k}$ be an intersecting family. 
If $r > 2^{2(1+H(\delta))}/\delta$ with  $\delta=\frac{1}{2\log_2(2k)}$, then $\hs^*(\hf,r)\not=\emptyset$. 
\end{fact}
\begin{proof}
 If $\hs(\hf,r)=\emptyset$, then $\hf$ is an $r$-spread family. Let $U_1$ be a $\frac{1}{2}$-random subset of $[n]$ and let $U_2=[n]\setminus U_1$.  Take $m=\log_2(2k)$, then we have $m\delta =\frac{1}{2}$. It follows from $r\delta > 2^{2(1+H(\delta))}$ that
\[
\left(\frac{1+H(\delta)}{\log_2(r\delta)}\right)^m < \frac{1}{2k}.
\]
By Theorem  \ref{thm-key} and the union bound, with positive probability that there exist $F_1, F_2 \in \hf$ such that $F_1 \subset U_1$, $F_2\subset U_2$. Then $F_1, F_2$ are disjoint, contradicting the fact that $\hf$ is intersecting. Thus $\hs(\hf,r)\neq \emptyset \neq \hs^*(\hf,r)$.
\end{proof}

We also need the following inequalities.

\begin{fact}\label{fact-2.2}
Let $n \ge 500$, $r= n/3$, $\delta=\frac{1}{2\log_2(2n)}$, $q= 4\log_2 n$ and $f(x) =r^{x}(n-x)!$, then

(i) $H(\delta) \le 0.288$ and $r \ge 2^{2(1+H(\delta))}/\delta$;

(ii) for any $q \le x\le n$, $f(x) \le (n-4)!$.
\end{fact}

\begin{proof}
Note that when $n\ge 500$, $\delta \leq \frac{1}{2\log_2(1000)}<\frac{1}{2}$, therefore, we have
 \[
 H(\delta) \leq H\left(\frac{1}{2\log_2(1000)}\right) \leq 0.288,
 \] 
and $r \ge 2^{2(1+H(\delta))}/\delta$.
 Since 
\[
\frac{f(i+1)}{f(i)} =\frac{n}{3(n-i)},
\]
we infer that $f(i+1)<f(i)$ for $i<\frac{2n}{3}$ and $f(i+1)>f(i)$ for $i\geq \frac{2n}{3}$. 
It follows that $f(x)\leq \max\{f(q),f(n)\}$ if $q \le x\le n$.
Note that 
\[
f(q)=  \left(\frac{n}{3}\right)^{q} (n-q)! < 2^{-q}\left(\frac{2n}{3}\right)^{q}(n-q)!<2^{-q}n!<n^{-4}n!<(n-4)!,
\]
since $2n/3 \le n-q$ for $n \ge 100$. Moreover for $n\geq 500$,
\[
f(n)= \left(\frac{n}{3}\right)^{n}  = \left(\frac{e}{3}\right)^n\left(\frac{n}{e}\right)^{n} \leq  n^{-4}n!\leq (n-4)!. 
\]
\end{proof}

For $\hb\subset \hp(X)$, define 
\[
\langle \hb \rangle = \{ F \in \hs_n \colon  \mbox{ there exists } B\in \hb \mbox{ such that }B\subset F \}.
\]

The following theorem was proved by Kupavskii and Zakharov in a more general setting. We slightly simplify the proof and improve some parameters. 

\begin{theorem}[\cite{kup}]\label{thm-key1}
Let $n \ge 500$ and $q=4\log_2 n$. For any intersecting family $\mathcal{F} \subset \mathcal{S}_n$, there exists an intersecting family $\mathcal{B}$ of proper subsets of size at most $q$ and $\mathcal{F}' \subset \mathcal{F}$ such that

(i) for any $B \in \mathcal{B}$ there is a family $\mathcal{F}_B \subset \mathcal{F}$ such that $\mathcal{F}_B(B)$ is $\frac{n}{3}$-spread;

(ii)  $\mathcal{F}\setminus \mathcal{F}' = \bigcup_{B\in \mathcal{B}}\mathcal{F}_B$;

(iii) $|\mathcal{F}'| \le (n-4)!$.
\end{theorem}

\begin{proof}
Let $r=n/3$, $\delta=\frac{1}{2\log_2(2n)}$ and $\hf_1=\hf$. In the first stage, let $\hb=\emptyset$.  Since $\mathcal{F}_1$ is intersecting and  $n\geq 500$,
by Facts \ref{fact-2.1} and \ref{fact-2.2}, we infer $\hs^*(\hf_1,r)\neq \emptyset$. 

If $|B|\leq q$ for all $B\in \hs^*(\hf_1,r)$, then  choose some $B_1\in \hs^*(\hf_1,r)$ and put it into $\hb$. Let 
\[
\hf_{B_1}=\{F\in \hf_1\colon B_1\subset F\}\mbox{ and } \hf_2=\hf_1\setminus \hf_{B_1}.
\]
Since $B_1$ is maximal and not $r$-spread in $\hf_1$, for any $C\supset B_1$, we have
\[
\frac{|\hf_{B_1}(C)|}{|\hf_{B_1}(B_1)|} = \frac{|\hf_{B_1}(C)|}{|\hf_1|}\times \frac{|\hf_1|}{|\hf_{B_1}(B_1)|} \le r^{-|C|+|B_1|}=r^{-(|C|-|B_1|)}.
\]
That is, $\hf_{B_1}(B_1)$ is $r$-spread. 
If $|B|\leq q$ for all $B\in \hs^*(\hf_2,r)$, then  choose some $B_2\in \hs^*(\hf_2,r)$ and put it into $\hb$. Let 
\[
\hf_{B_2}=\{F\in \hf_2\colon B_2\subset F\}\mbox{ and } \hf_3=\hf_2\setminus \hf_{B_2}.
\]
Similarly $\hf_{B_2}(B_2)$ is $r$-spread.

Do it repeatedly until there exists some $B_i\in \hs^*(\hf_i,r)$ with $|B_i|\geq q+1$ or  
$\hf_i=\emptyset$. Suppose that the procedure ends at $\hf_\ell$ and  let $\hf'=\hf_\ell$.  
If  there exists some $B_\ell\in \hs^*(\hf_\ell,r)$ with $|B_\ell|\geq q+1$, then
\[
|\hf_{\ell}(B_\ell)|>r^{-|B_\ell|} |\hf_{\ell}|.
\]
Let $|B_\ell|=x$. Clearly $q< x\leq n$.
\[
|\hf'|=|\hf_\ell| \leq r^x |\hf_{\ell}(B_\ell)| \leq  \left(\frac{n}{3}\right)^{x} (n-x)!=:f(x). 
\]
Then by Fact \ref{fact-2.2}, we have $|\hf'|\leq (n-4)!$ and (ii), (iii) hold.

Note that by Fact \ref{fact-2.1}, for any $B_i \in \hb$, $B_i \not=\emptyset$.  Now we show that $\hb$ is intersecting. Suppose not, then there exists $B_i,B_j\in \hb$ such that $B_i\cap B_j=\emptyset$. 
Recall that $\hf_i(B_i)$, $\hf_j(B_j)$ are both $r$-spread. Let $R=B_i\cup B_j$ and let 
\[
\hg_i=\left\{G_i \in \hf_i(B_i) \colon G_i \cap B_j = \emptyset  \right\}
\]
and
\[
\hg_j=\left\{G_j \in \hf_j(B_j) \colon G_j \cap B_i = \emptyset  \right\}.
\]
For $n\geq 500$, 
\[
|\hg_i| \geq |\hf_i(B_i)| -\sum_{x\in B_j} |\hf_i(B_i\cup \{x\})| \geq |\hf_i(B_i)| - \frac{q}{r} |\hf_i(B_i)| \geq \left(1-\frac{4\log_2 n}{n/3}\right) |\hf_i(B_i)|\geq \frac{3|\hf_i(B_i)|}{4}. 
\]
For any $S\in \hp(X\setminus R)$, we have
\[
|\hg_i(S)| \leq |\hf_i(B_i\cup S)| \leq r^{-|S|} |\hf_i(B_i )|\leq r^{-|S|}  \frac{4}{3}|\hg_i|  \leq \left(\frac{3r}{4}\right)^{-|S|} |\hg_i|.
\]
It follows that $\hg_i$ is $\frac{3r}{4}$-spread. Similarly, $\hg_j$ is also $\frac{3r}{4}$-spread.
Let $Y_1$ be a $\frac{1}{2}$-random subset of $X\setminus R$ and let $Y_2=(X\setminus R)\setminus Y_1$.  Let $m=\log_2(2n)$ and recall $\delta=\frac{1}{2\log_2(2n)}$. Then $m\delta =\frac{1}{2}$. Since $n\geq 500$, $\frac{3r}{4}\delta = \frac{n}{8\log_2(2n)}\geq 6$. Note that $H(\delta) \leq 0.288$ by Fact \ref{fact-2.2}.
Using Theorem  \ref{thm-key}, we have 
 \[
 Pr[\exists G_i\subset Y_1 \mbox{ for some }G_i\in \hg_i]>  1- \left(\frac{1+H(\delta)}{\log_2(3r\delta/4)}\right)^m n>\frac{1}{2}. 
 \]
 and similarly,
 \[
 Pr[\exists G_j\subset Y_2 \mbox{ for some }G_j\in \hg_j]>  1- \left(\frac{1+H(\delta)}{\log_2(3r\delta/4)}\right)^m n>\frac{1}{2}. 
 \]

 It follows that with positive probability that there exist $G_i\in \hg_i$, $G_j\in \hg_j$ such that $G_i\subset Y_1$ and $G_j\subset Y_2$. Then $G_i\cup B_i, G_j\cup B_j\in \hf$ are disjoint, contradicting the fact that $\hf$ is intersecting. Thus $\hb$ is intersecting. 
\end{proof}

\section{Proof of Theorem \ref{main_thm}}

We say $F_0,F_1,F_2,\ldots,F_{s}\subset X$ form a  {\it pseudo sunflower of size $s+1$ with center $C$} if $C\subsetneq F_0$ and $F_0\setminus C,F_1\setminus C,\ldots, F_{s}\setminus C$ are pairwise disjoint. 
The following result on pseudo sunflowers was proved by F\"{u}redi \cite{Fu80}.

\begin{theorem}[F\"{u}redi \cite{Fu80}, cf. also \cite{frankl}]\label{thm-furedi}
Let $\mathcal{F}$ be a $k$-uniform family without any pseudo sunflower of size $s+1$, then
\[
|\mathcal{F}| \le s^k.
\]
\end{theorem}

Let us mention that the term pseudo sunflower originates from \cite{frankl}.

For a family $\hh\subset \hp(X)$ and $s \in [n]$,  we define a basis $\ha(\hh,s)$ of $\hh$ by a pseudo sunflower procedure introduced by Kupavskii and Zakharov \cite{kup}. We start with $\hh_0=\hh$. For $i=0,1,2,\ldots$, if $\hh_i$ contains an $s$-uniform pesudo sunflower of size $s+1$ with center $C$, then let $\hh_{i+1}$ be a family obtained from $\hh_i$ by removing all sets containing $C$ and  adding the center $C$. We repeat the procedure until the resulting family does not  contain any pesudo sunflower of size $s+1$. Denoting the resulting family as $\ha(\hh,s)$.

\begin{claim}
If $\hh$ is an intersecting family of subsets with size at most $s$, then $\ha(\hh,s)$ is also intersecting. 
\end{claim}
\begin{proof}
Suppose not. Let $S_1,S_2\in \ha(\hh,s)$ and $S_1\cap S_2 = \emptyset$.
Then at least one of $S_1,S_2$ is the center of some pseudo sunflower in $\hh$. Otherwise  $\hh$ is not intersecting. Say $S_1$ is the center of pseudo sunflower $F_0,F_1,F_2,\ldots,F_{s}\subset \hh$ such that $S_1\subsetneq F_0$ and $F_0\setminus S_1,F_1\setminus S_1,\ldots, F_{s}\setminus S_1$ are pairwise disjoint. 
Since $S_1\cap S_2 = \emptyset$, for $i=0,1,\ldots, s$, $(F_i\setminus S_1)\cap S_2 \not= \emptyset$. However, $|S_2| \le s$, a contradiction.
\end{proof}

\noindent{\bf Proof of Theorem \ref{main_thm}:~~} Let $\hf\subset \hs_n$ be an intersecting family. By Theorem \ref{thm-key1}, there exists an intersecting family $\hb\subset \hp(X,q)$ such that (i),(ii),(iii) of Theorem \ref{thm-key1} hold. 

Let $\hht$ be the minimal members in $\hb$. Let 
$\hp_q=\ha(\hht,q)$ and $\ha_{q}$ be the set of all $q$-subsets in $\hp_q$.  Then let $\hp_{q-1} =\ha(\hp_q\setminus \ha_q,q-1)$ and  $\ha_{q-1}$ be the set of all $(q-1)$-subsets in $\hp_{q-1}$. Define $\ha_{q-2},\ldots,\ha_4,\ha_3$ and let $\ha_2= \hp_3\setminus \ha_3$.  Clearly $\ha_2$ is intersecting. By Theorem \ref{thm-furedi} we obtain that
\[
|\ha_i|\leq i^i,\ i=3,4,\ldots,q.
\]
Since $\frac{(i+1)^{i+1}(n-i-1)!}{i^i(n-i)!} = (1+\frac{1}{i})^i \frac{i+1}{n-i}\leq e\frac{q+1}{n-q} \leq \frac{1}{2}$ for $n\geq 500$, we have
\[
\sum_{3\leq i\leq q} |\ha_i|(n-i)! \leq 2\times 3^3(n-3)! =54 (n-3)!.
\]

If $\ha_2$ is a star of center $x$, then for $n\geq 100$, 
\[
\gamma(\mathcal{F} )\leq |\hf(\bar{x})|\le \sum_{3\leq i\leq q} |\ha_i|(n-i)! +|\mathcal{F}'| \leq  54 (n-3)!+(n-4)!<(n-3)(n-3)!. 
\]
Thus we may assume $\ha_2$ is non-trivial intersecting.  It follows that $\ha_2$ is a triangle. Without loss of generality assume that 
\[
\ha_2=\{\{x_{11},x_{22}\},\{x_{11},x_{33}\},\{x_{22},x_{33}\}\}. 
\]
Let $T=\{x_{11},x_{22},x_{33}\}$. 
If $|F\cap T|\geq 2$ for all $F\in \hf$, then $\gamma(\mathcal{F} )\leq (n-3)(n-3)!$ follows. 
Thus we may assume that there exists $F_0\in \hf$ such that $|F_0\cap T|\leq 1$. By symmetry assume that $\{x_{22},x_{33}\}\cap F_0=\emptyset$. 
Note that each $F$ containing $\{x_{22},x_{33}\}$ has to  intersect $F_0$. 
Clearly, the number of $F$ containing $\{x_{22},x_{33}\}$ is $(n-2)!$, and the number of such $F$ does not intersect $F_0$ is at least the number of derangements of $S_{n-2}$, i.e., the integer part of $((n-2)!+1)/e$. It follows that
\[
| \langle \{x_{22},x_{33}\} \rangle|\leq (1-\frac{1}{e})(n-2)!
\]
Then, when $n \geq 500$, $(1-1/e)(n -2) + 55 < n-3$, therefore,
\begin{align*}
\gamma(\mathcal{F})\leq |\hf(\overline{x_{11}})|&\le (1-\frac{1}{e})(n-2)!+ \sum_{i = 3}^{q} |\ha^{(i)}|(n-i)! + |\mathcal{F}'|\\[3pt]
&\leq (1-\frac{1}{e})(n-2)!+ 54 (n-3)!+(n-4)! \\[3pt]
&<(n-3)(n-3)!.
\end{align*}

This completes the proof of Theorem \ref{main_thm}. $\hfill \Box$

\noindent{\bf Acknowledgments.} 
We would like to thank the anonymous referees for their helpful comments. The first author is supported by National Natural Science Foundation of China (No. 12471316) and  Natural Science Foundation of Shanxi Province, China (Nos. RD2500002993, RD2200004810). 
The second author is supported by National Natural Science Foundation of China (Nos. 12501460, 12131011) and Shenzhen Science and Technology Program (No. RCBS20221008093102011).

\end{document}